\documentclass[10pt,reqno]{amsart}
\usepackage{amsmath}
\usepackage{amsthm}
\usepackage{amssymb}
\usepackage{amsfonts}
\usepackage{epsfig}
\usepackage{epsf}
\newtheorem{theorem}{Theorem}

\newtheorem{example}{Example}
\newtheorem*{remark}{Remark}
\newtheorem{lemma}{Lemma}

\newcommand{\baseRing}[1]{\ensuremath{\mathbb{#1}}}
\newcommand{\C}{\baseRing{C}}

\newcommand{\N}{\baseRing{N}}
\newcommand{\g}{\mathfrak{g}}
\newcommand{\h}{\mathfrak{h}}

\renewcommand{\phi}{\varphi}

\title{\bf Constructing $r-$matrices on Simple Lie Superalgebras}
\author{Gizem Karaali}
\address{Department of Mathematics, University of California, Berkeley, 
Berkeley, CA 94720, USA}
\email{gizem@math.berkeley.edu}

\begin{document}

\abstract{We construct $r-$matrices for simple Lie superalgebras
with non-degenerate Killing forms using Belavin-Drinfeld type
triples. This construction gives us the standard $r-$matrices and
some nonstandard ones.}

\endabstract

\maketitle

\section{Introduction}
\label{SectionIntroduction}

Let $\g$ be a Lie algebra with a non-degenerate
$\g-$invariant bilinear form ${(  \; , \; )}$. Denote by
$\Omega$ the element of ${(\g \otimes \g)^{\g}}$ that corresponds
to the quadratic Casimir element in $\mathfrak{Ug}$ of $\g$.
Then the \textbf{classical Yang-Baxter equation} (\textbf{CYBE})
for an element ${r \in \g \otimes \g}$ is:  
$$[r^{12},r^{13}] + [r^{12},r^{23}] + [r^{13},r^{23}] = 0. $$ 
\noindent 
A solution $r$ to the classical Yang-Baxter equation
is called a \textbf{classical $r-$matrix} (or simply an
\textbf{$r-$matrix}). An  $r-$matrix $r$ is called
\textbf{non-degenerate} if it satisfies: 
$$r^{12} + r^{21} = \Omega.$$
\noindent
In \cite{BD1} and \cite{BD2} Belavin and Drinfeld classified such 
$r-$matrices. The solutions in this classification are parametrized by 
triples ${(\Gamma_1, \Gamma_2, \tau)}$ (called \textbf{admissible 
triples}) where the $\Gamma_i$ are certain
subsets of the set $\Gamma$ of simple roots of $\g$, and
${\tau : \Gamma_1 \rightarrow \Gamma_2}$ is an isometric
bijection. They proved  that for each admissible triple and some
fixed ${r_0 \in \g \otimes \g}$  there exists a unique
non-degenerate classical $r-$matrix, and  conversely that each
non-degenerate classical $r-$matrix can be  associated with such
data.

In this paper we construct classical $r-$matrices
using analogs of the Belavin-Drinfeld data for simple Lie
superalgebras with non-degenerate Killing form. We first start
with a review of the situation in the Lie algebra case. Thus in
Section \ref{SectionLieAlgebraCase} we give an overview of the
Belavin-Drinfeld result for simple Lie algebras. Next, in
Section \ref{SectionSuperTheory}, we recall some basic
definitions and results about simple Lie superalgebras. Then
after developing the necessary ingredients we state our main
theorem. The next three sections of the paper are devoted to the
proof of this theorem. Then in Section \ref{Sectiononsl21} we
construct various $r-$matrices for the Lie superalgebra
${sl(2,1)}$ using the main theorem. 

This theorem is very much in the spirit of the Belavin-Drinfeld
result. It tells us that, given a Belavin-Drinfeld type triple,
one can construct a non-degenerate $r-$matrix in a way similar to
the construction in the Lie algebra case. One should point out,
however, that this is not a complete classification result. The
theorem gives us $r-$matrices, but does not tell us whether we
can always extract a Belavin-Drinfeld type triple from a given
$r-$matrix. We discuss this briefly in the last section of the paper. 

\subsection*{Acknowledgments.} I would like to thank Nicolai
Reshetikhin for fruitful discussions, Vera Serganova for her
close attention and  useful suggestions, and Milen Yakimov for
introducing me to the  Belavin-Drinfeld result. This work was
supported by the NSF grant  DMS$-0070931$.

\section{Classification Theorem for Lie Algebras}
\label{SectionLieAlgebraCase}

Here we recall briefly the main result of \cite{BD1} and \cite{BD2}
for Lie algebras.
Let $\g$ be a simple Lie algebra. Fixing a positive Borel
subalgebra $\mathfrak{b}_{+}$ determines a Cartan subalgebra
$\h$ for $\g$, and we can talk about positive roots, simple
roots etc. Hence we can define ${\Gamma = \{ \alpha_1, \alpha_2,
\cdots, \alpha_r \}}$ to be the set of all simple roots of
$\g$. We will be interested in \textbf{admissible
triples}, i.e. triples $(\Gamma_1, \Gamma_2, \tau)$ where
$\Gamma_i \subset \Gamma$ and ${\tau : \Gamma_1 \rightarrow
\Gamma_2}$ is a bijection such that 

\begin{enumerate} 

\item for any $\alpha, \beta \in \Gamma_1$,
$(\tau(\alpha),\tau(\beta)) = (\alpha,\beta)$;

\item for any $\alpha \in \Gamma_1$ there exists a $k \in \N$
such that $\tau^k(\alpha) \not \in \Gamma_1$. \footnotemark
\footnotetext{The expression $\tau^k(\alpha)$ has a meaning only
if the expressions $\tau^j(\alpha)$ for all $j < k$ are
elements of $\Gamma_1$. So this condition actually may be
translated as saying that $\tau^k(\alpha)$ does not make
sense for large enough $k$.} 

\end{enumerate}

We will also need a continuous parameter $r_0$, an element of
$\h \otimes \h$ which satisfies the following equations: 
\begin{equation}
\label{unitarityinh}
{r_0}^{12} + {r_0}^{21} = \Omega_0 
\end{equation}
\begin{equation}
\label{YBinh}
(\tau(\alpha) \otimes 1) (r_0) + (1 \otimes \alpha) (r_0) = 0
\textmd{  for  all } \alpha \in \Gamma_1 
\end{equation}

\noindent
where $\Omega_0 \in \h \otimes \h$ is the $\h-$component of
the quadratic Casimir element $\Omega$ of $\g$. 

Fix a system of
Weyl-Chevalley generators $X_{\alpha}, Y_{\alpha},
H_{\alpha}$ for $\alpha \in \Gamma$. Recall that these elements 
generate the Lie algebra $\g$ with the defining relations: 
$[X_{\alpha_i}, Y_{\alpha_j}] =
\delta_{i,j} H_{\alpha_j}$, $[H_{\alpha_i}, X_{\alpha_j}] =
a_{i,j} X_{\alpha_j}$ and $[H_{\alpha_i}, Y_{\alpha_j}] = -
a_{i,j} Y_{\alpha_j}$ for all $\alpha_i, \alpha_j \in \Gamma$,
where $a_{i,j} = \alpha_j(H_{\alpha_i}) = \frac{2 (\alpha_i,
\alpha_j)}{(\alpha_i,\alpha_i)}$, along with the well-known Serre 
relations. 

Denote by $\g_i$ the subalgebra of $\g$ generated by the
elements $X_{\alpha}, Y_{\alpha}, H_{\alpha}$ for all $\alpha
\in \Gamma_i$. We define a map $\varphi$ by: 

$$ \varphi(X_{\alpha}) =
X_{\tau(\alpha)}  \; \; \; \; \; \; \; \varphi(Y_{\alpha}) =
Y_{\tau(\alpha)}  \; \; \; \; \; \; \;  \varphi(H_{\alpha}) =
H_{\tau(\alpha)} $$ 
\noindent 
for all $\alpha \in \Gamma_1$. Then this can be extended to an
isomorphism ${\varphi : \g_1 \rightarrow \g_2}$ because the
relations between $X_{\alpha}, Y_{\alpha}, H_{\alpha}$ for
$\alpha \in \Gamma_1$ will be the same as the relations between 
$X_{\tau(\alpha)}, Y_{\tau(\alpha)}, H_{\tau(\alpha)}$ for
$\alpha \in \Gamma_1$. Note that this requires the first
property of $\tau$, namely that it is an isometry. Next 
extend $\tau$ to a bijection between the $\overline{\Gamma}_i$,
where $\overline{\Gamma}_i$ is the set of those roots which can
be written as a nonnegative integral linear combination of the
elements of $\Gamma_i$. In each root space $\g_{\alpha}$, we
choose an element $e_{\alpha}$ such that
$(e_{\alpha},e_{-\alpha}) = 1$ for any $\alpha$ and
$\varphi(e_{\alpha}) = e_{\overline{\tau}(\alpha)}$ for all
$\alpha \in \overline{\Gamma}_1$ \footnotemark.  
\footnotetext{ Such $e_{\alpha}$ can always be chosen
consistently if ``there are no cycles ", i.e. if $\tau$ satisfies
the second property. Otherwise, if there is a cycle of simple
roots $\alpha_1, \cdots, \alpha_k$ such that $\tau(\alpha_i) =
\alpha_{i+1}$ for $i = 1, 2, \cdots, k-1$ and $\tau(\alpha_k) =
\alpha_1$ and $\tau$ is an isometry, then we must have that
$(\alpha_i, \alpha_{i+1}) = 0$ for all $i = 1,2, \cdots, k-1$,
since no cycles are allowed in Dynkin diagrams. Then if $s$ is
the smallest integer such that $(\alpha_1, \alpha_{1+s}) \ne 0$,
then $s > 1$ and $s | k$. So our cycle must have $s$
disconnected subgraphs of length $k/s$. Then we can still choose
$e_{\alpha}$  consistently, provided we allow only cycles as
above. }
Finally we define a partial order on the set of all positive
roots:    
$$ \alpha < \beta \textmd{ if and only if there exists a } k
\in \N \textmd{ such that } \beta = \overline{\tau}^k(\alpha)$$    
\noindent 
Note that in this case we will have $\alpha \in
\overline{\Gamma}_1, \beta \in \overline{\Gamma}_2$. Clearly
one needs property $2$ for a partial order; no cycles are
allowed in partial orders. 
\\
\\
Now we can state the Belavin-Drinfeld theorem:

\begin{theorem}
\label{LieTheorem}
(1) Let $r_0 \in \h \otimes \h$ satisfy Equations \ref{unitarityinh} and 
\ref{YBinh}. Then the element $r$ of $\g \otimes \g$ defined by:
$$ r = r_0 + \sum_{\alpha > 0} e_{-\alpha} \otimes e_{\alpha} 
+ \sum_{\alpha,\beta > 0, \alpha < \beta} (e_{-\alpha} \otimes
e_{\beta} - e_{\beta} \otimes e_{-\alpha}) $$
\noindent
is a solution to the system:
\begin{equation*}
\tag{\ref{unitarity}}
r^{12} + r^{21} = \Omega 
\end{equation*}
\begin{equation*}
\tag{\ref{YBequation}}
[r^{12},r^{13}] + [r^{12},r^{23}] + [r^{13},r^{23}] = 0
\end{equation*}
(2) Any solution (up to isomorphism) to the above system can be
obtained as above by a suitable choice of an admissible
triple $(\Gamma_1,\Gamma_2,\tau)$ and some $r_0 \in \h \otimes
\h$ that satisfies Equations \ref{unitarityinh} and \ref{YBinh}.  
\end{theorem}

The proof of this theorem provided in \cite{BD2} is actually quite
clear; one can also look at \cite{ES} for another exposition. 

\section{The Construction Theorem for Lie Superalgebras}
\label{SectionSuperTheory}

Now our aim is to develop a similar theory for super
structures. Let $\g$ be a simple Lie superalgebra with a
non-degenerate Killing form \footnotemark.
\footnotetext{ This implies that $\g$ is isomorphic to a simple
Lie algebra or to one of the following classical Lie
superalgebras: $A(m,n)$ with $m \ne n$, $B(m,n)$, $C(n)$,
$D(m,n)$ with $m - n \ne 1$, $F(4)$, or $G(3)$. See \cite{K} for
details.} 

\subsection{The Quadratic Casimir Element:} Let $\{ I_{\alpha}
\}$ be a homogeneous basis for $\g$ and denote by $\{
{I_{\alpha}}^* \}$ the dual basis of $\g$ with respect to the
non-degenerate (Killing) form. Thus we have: 
$$ (I_{\alpha},{I_{\beta}}^*) = \delta_{\alpha \beta}$$
\noindent
If we denote the parity of a homogeneous element $x \in \g$ by $|x|$, 
then we also have that 
$$ |I_{\alpha}| = |{I_{\alpha}}^*| $$
\noindent
since the supertrace form is consistent [Recall that a bilinear
form $($ , $)$ is \textbf{consistent} if for any homogeneous $x,y
\in \g$ of different parities, $(x,y) = 0$]. 

We can write the quadratic Casimir element of $\g$ as follows:
$$ \Omega = \sum_{\alpha} (-1)^{|I_{\alpha}||{I_{\alpha}}^*|}
I_{\alpha} \otimes {I_{\alpha}}^*  = \sum_{\alpha}
(-1)^{|I_{\alpha}|} I_{\alpha} \otimes {I_{\alpha}}^* $$

\begin{example}
\label{Exampleglmn}
For the special case $\g = gl(m,n)$ with basis {${\{
e_{i,j} | 1 \le i,j \le m+n\}}$}, we use the supertrace form to
find the dual basis:  $$ {e_{i,j}}^* = (-1)^{[i]} e_{j,i}$$
\noindent
where
$$ [j] = \left \{ \begin{matrix} 0 & if & j \le m \\1 & if & j > m
\end{matrix} \right .$$
\noindent
Then this gives us:
$$ \Omega = \sum_{\alpha} (-1)^{|I_{\alpha}|} I_{\alpha} \otimes
{I_{\alpha}}^* = \sum_{i,j} (-1)^{|e_{i,j}|} 
e_{i,j} \otimes (-1)^{[i]} e_{j,i} = \sum_{i,j} (-1)^{[j]}
e_{i,j} \otimes e_{j,i} $$
\end{example}

\subsection{Borel subsuperalgebras and Dynkin diagrams: } Let $\h$ 
be a Cartan subalgebra. By definition, $\h \subset
\g_{\overline{0}}$ is a Cartan subalgebra of the even
part of our Lie superalgebra. Let $\Delta =
\Delta_{\overline{0}} + \Delta_{\overline{1}}$ be the set of
all roots of $\g$ associated with the Cartan subalgebra $\h$ 
\footnotemark, where $\Delta_{\overline{0}}$ and
$\Delta_{\overline{1}}$ are the even and odd roots
respectively. Fix a Borel subsuperalgebra
$\mathfrak{b}$ containing $\h$. Recall that a Lie
subsuperalgebra $\mathfrak{b}$ of a Lie superalgebra $\g$ is a
\textbf{Borel subsuperalgebra} if there is some Cartan
subsuperalgebra $\h$ of $\g$ and some base $\Gamma$ for
$\Delta$, such that   
$$ \mathfrak{b} = \h \oplus \bigoplus_{\alpha \in \Delta^+}
\g_{\alpha} $$   
\noindent 
where $\Delta^+$ consists of all nonnegative integral
combinations of the elements of $\Gamma$ that are in $\Delta$.

\footnotetext{In fact $\Delta$ is independent of the choice of $\h$}

In the Lie algebra case subalgebras given by this definition
turn out to be maximally solvable, and all maximally
solvable subalgebras of a simple Lie algebra are of this type.
Thus this definition agrees with the usual definition of a
Borel subalgebra as a maximally solvable subalgebra. However
Borel subsuperalgebras as defined above are not necessarily
maximally solvable.

\begin{example}
Let $\g$ be a Lie superalgebra, fix some set $\Delta_+$ of
positive roots, and let $\alpha$ be a positive isotropic root. 
Define $\mathfrak{b}$ as the sum of all the positive
root spaces. Then $\mathfrak{b}$ is a Borel subsuperalgebra but
is not maximally solvable. The (parabolic \footnotemark)
subsuperalgebra ${\mathfrak{p} = \mathfrak{b} \oplus
\g_{-\alpha}}$ is also solvable.  
\end{example}

\footnotetext{As in the Lie algebra case, a Lie subsuperalgebra
$\mathfrak{p}$ is a \textbf{parabolic subsuperalgebra} if
$\mathfrak{p}$ contains a Borel. } 

In fact maximally solvable subsuperalgebras may be even more
complicated. [See \cite{Shch} for a study of maximally solvable
subsuperalgebras of ${gl(m,n)}$ and ${sl(m,n)}$.] Therefore we
choose to define Borels as above instead of using the more
traditional characterization by maximal solvability. 

Thus Borel subsuperalgebras determine simple roots, and different Borel 
subsuperalgebras may correspond to different Dynkin diagrams and Cartan 
matrices. Let us then fix some Borel subsuperalgebra ${\mathfrak{b}}$, or 
equivalently some set $\Gamma$ of simple roots, and the
associated Dynkin diagram $D$. Note that $\Gamma$ may contain
even and odd roots. Another significant difference from the
theory of Lie algebras is to be noted here; two Dynkin diagrams
of a given Lie superalgebra are not necessarily isomorphic, but
can be obtained from one another via a chain of odd reflections. 

\subsection{The Data for the Theorem: } In this setup, let
$\Delta^+$ (resp. $\Delta^-$) be the set of all positive (resp.
negative) roots, with respect to the chosen ${\Gamma = \{
\alpha_1, \alpha_2, \cdots, \alpha_r  \}}$. 

Now let $\Gamma_1, \Gamma_2 \subset \Gamma$ be two subsets, and
${\tau : \Gamma_1 \rightarrow \Gamma_2}$ be a bijection. We will
say that the triple ${(\Gamma_1,\Gamma_2,\tau)}$ is  
\textbf{admissible} if:

\begin{enumerate}

\item for any $\alpha, \beta \in \Gamma_1$,
$(\tau(\alpha),\tau(\beta)) = (\alpha,\beta)$;

\item for any $\alpha \in \Gamma_1$ there exists a $k \in \N$
such that $\tau^k(\alpha) \not \in \Gamma_1$;

\item $\tau$ preserves grading, i.e. it maps even roots to even
ones, and odd roots to odd ones.

\end{enumerate}

For a given admissible triple $(\Gamma_1, \Gamma_2, \tau)$, we
define $\overline{\Gamma}_i$ for $i=1,2$ as in the Lie algebra
case: $\overline{\Gamma}_i$ is the set of those roots which can
be written as a nonnegative integral linear combination of the
elements of $\Gamma_i$. Then we can extend $\tau$ linearly to a
bijection $\overline{\tau} : \overline{\Gamma}_1 \rightarrow
\overline{\Gamma}_2$. Using $\overline{\tau}$, we define a
partial order on $\Delta^+$:  
$$ \alpha < \beta \textmd{ if and only if there exists a } k \in
\N \textmd{  such that } \beta = \overline{\tau}^k(\alpha)$$  

For any $\alpha \in \Gamma$, pick a nonzero $e_{\alpha} \in
\g_{\alpha}$. Then since each of the $\g_{\alpha}$ are
one dimensional, and the Killing form is a non-degenerate
pairing of $\g_{\alpha}$ with {$\g_{-\alpha}$\footnotemark,} we
can uniquely pick $e_{-\alpha} \in \g_{-\alpha}$ such that: 
$$(e_{\alpha}, e_{-\alpha}) = 1. $$

\footnotetext{ This will hold for all classical Lie
superalgebras of the form: $A(m,n)$ for $(m,n) \ne (1,1)$,
$B(m,n)$, $C(n)$, $D(m,n)$, $D(2,1; \alpha)$, $F(4$, and $G(3)$.
In fact the dimension of $\g_{\alpha}$ is one provided $\g$ is a
classical Lie superalgebra different from $A(1,1)$, $P(2)$,
$P(3)$, and $Q(n)$. See \cite{K} for more on Lie superalgebras.}

Then we will have, for each $\alpha \in \Gamma$:
$$ [e_{\alpha}, e_{-\alpha}] = (e_{\alpha}, e_{-\alpha})
h_{\alpha}$$
\noindent
where $h_{\alpha}$ is the nonzero vector
defined by $(h_{\alpha},h) = \alpha(h)$ for all $h \in \h$.
The set ${\{ h_{\alpha} | \alpha \in \Gamma \}}$ is a basis for
$\h$. Hence we can write $\Omega_0$, the $\h-$part of the
quadratic Casimir $\Omega$ of $\g$, as follows:
$$ \Omega_0 = \sum_{i=1}^{r} h_{\alpha_i} \otimes
{h_{\alpha_i}}^* $$
\noindent
where the set $\{ {h_{\alpha}}^* | \alpha \in \Gamma \}$ is the
basis in $\h$ dual to ${\{ h_{\alpha} | \alpha \in \Gamma \}}$

Actually we can choose a nonzero $e_{\alpha} \in \g_{\alpha}$
for each $\alpha \in \Delta$ such that ${(e_{\alpha},
e_{-\alpha}) = 1}$ whenever $\alpha$ is positive. Then we compute
the duals with respect to the standard (Killing) form:   
$${e_{\alpha}}^* = e_{-\alpha}$$ 
$${e_{-\alpha}}^* =  (-1)^{|\alpha|} e_{\alpha}$$    
\noindent
for all positive roots $\alpha$. Here $|\alpha|$ is the parity
of the root $\alpha$.

Then we can see that the quadratic Casimir
element of our Lie superalgebra $\g$ will be: 
\begin{eqnarray*} 
\Omega &=& \sum_i (-1)^{|I_i|} I_i \otimes {I_i}^*  \\
\\
&=& \sum_{i=1}^{r} h_{\alpha_i} \otimes {h_{\alpha_i}}^*
+ \sum_{\alpha \in \Delta} (-1)^{|e_{\alpha}|} e_{\alpha}
\otimes {e_{\alpha}}^* \\
\\
&=& \Omega_0 + 
\sum_{\alpha \in \Delta^+} (-1)^{|\alpha|} e_{\alpha} \otimes
e_{-\alpha} + 
\sum_{\alpha \in \Delta^+} e_{-\alpha} \otimes e_{\alpha}  
\end{eqnarray*}

\textbf{Example \ref{Exampleglmn} continued: }
\textit{Let us consider the special case $\g = gl(m,n)$ again.
We can see that the positive root spaces will correspond to
$e_{i.j}$ for $i < j$. So if we choose $e_{\alpha}$s for the
positive root $\alpha$ to be the $e_{i,j}$ in $\g_{\alpha}$, we
will have $i < j$ and $e_{-\alpha}$ will be $(-1)^{[i]}
e_{j,i}$. Then we will have:   
$$ {e_{\alpha}}^* = {e_{i.j}}^* = (-1)^{[i]}e_{j,i}
= e_{-\alpha} $$ 
$$ {e_{-\alpha}}^* = (-1)^{[i]} {e_{j,i}}^* =
(-1)^{[i]} (-1)^{[j]} e_{i,j} = (-1)^{|\alpha|} e_{\alpha} $$
\noindent
and the above formula for $\Omega$ will agree with the Casimir
element found earlier.}

\subsection{Statement of the Theorem: } We are now ready to state
our main theorem. Its proof will be presented in the next three
sections.

\begin{theorem}
\label{SuperTheorem}
Let $r_0 \in \h \otimes \h$ satisfy:
\begin{equation*}
\tag{\ref{unitarityinh}}
{r_0}^{12} + {r_0}^{21} = \Omega_0 
\end{equation*}
\begin{equation*}
\tag{\ref{YBinh}}
(\tau(\alpha) \otimes 1) (r_0) + (1 \otimes \alpha) (r_0) = 0
\textmd{  for  all } \alpha \in \Gamma_1 
\end{equation*}
\noindent
Then the element $r$ of $\g \otimes \g$ defined by:
\begin{equation*}
\tag*{$(*)$}
r = r_0 + 
\sum_{\alpha > 0} e_{-\alpha} \otimes e_{\alpha}  + 
\sum_{\alpha,\beta > 0, \alpha < \beta} 
(e_{-\alpha} \otimes e_{\beta} - (-1)^{|\alpha|} e_{\beta}
\otimes e_{-\alpha})  
\end{equation*} 
\noindent
is a solution to the system:
\begin{equation*}
\tag{\ref{unitarity}}
r^{12} + r^{21} = \Omega 
\end{equation*}
\begin{equation*}
\tag{\ref{YBequation}}
[r^{12},r^{13}] + [r^{12},r^{23}] + [r^{13},r^{23}] = 0
\end{equation*}
\end{theorem}

\begin{remark}
Comparing this theorem with the corresponding theorem for Lie
algebras, we notice that there is a sign introduced in the
equation defining the $r-$matrix associated to a given
admissible triple. However since Lie algebras are Lie
superalgebras with only even roots, Equation $(*)$
reduces to the corresponding equation in Theorem
\ref{LieTheorem} when the Lie superalgebra $\g$ in question is
actually a Lie algebra.     
\end{remark}

\section{Technical Lemmas}

Let $\g$ be a simple Lie superalgebra with a non-degenerate
Killing form \footnotemark. Assume that we fix a homogeneous
basis $\{ I_{\alpha} \}$ for $\g$ and denote by $\{
{I_{\alpha}}^* \}$ the dual basis of $\g$ with respect to the
non-degenerate (Killing) form. We will first prove:
\footnotetext{Our results can mostly be extended to
the class of classical Lie superalgebras. If $\g$ is classical,
there is a non-degenerate invariant form on $\g$. In the
following, most of the statements involving the Killing form may
be asserted more generally for such an invariant form.} 

\begin{lemma} 
\label{rstatementsintofstatements}
Let $r \in \g \otimes \g$ be such that
$$r = (f \otimes 1) \Omega $$
\noindent
Then the system of equations:
\begin{equation}
r^{12} + r^{21} = \Omega
\label{unitarity}
\end{equation}
\begin{equation}
[r^{12},r^{13}] + [r^{12},r^{23}] + [r^{13},r^{23}] = 0
\label{YBequation}
\end{equation}
\noindent
is equivalent to the system:
\begin{equation}
f + f^* = 1
\label{unitarityforf}
\end{equation}
\begin{equation}
(f-1)[f(x),f(y)] = f([(f-1)(x),(f-1)(y)])
\label{YBforf}
\end{equation}
\noindent
where $f^*$ stands for the adjoint of $f$ with respect to the standard
from $($ , $)$.
\end{lemma}

\begin{remark}
This is exactly the same technical lemma used in the proof of
the main classification theorem in the Lie algebra case as
presented in \cite{ES}. Our proof here will be a
generalization of the proof provided there. We will use general
properties like the invariance, non-degeneracy and the
consistency of the Killing form. The main difference will be
that in our case, we may not be able to pick an orthonormal
basis for our Lie superalgebra $\g$. However it turns out that a
pair of dual bases will be sufficient for our purposes. 
\end{remark}

\noindent
\textbf{Proof: }
\underline{PART 1}: [This part is exactly the same as in
the Lie algebra case.] We have:  
$$ r^{12} + r^{21} = (f \otimes 1)\Omega + (1 \otimes f)\Omega
= (f \otimes 1)\Omega + (f^* \otimes 1)\Omega = ((f +
f^*)\otimes 1)  \Omega. $$ 
\noindent
Then we have: 
$$ \begin{matrix}
\Omega = r^{12} + r^{21}\\
\Leftrightarrow \\
\Omega = (1 \otimes 1) \Omega = ((f + f^*)\otimes 1) \Omega \\
\Leftrightarrow \\
 1 = (f + f^*)
\end{matrix}$$
\noindent
This proves the equivalence of the statements: 
$$ \begin{matrix} 
\Omega = r^{12} + r^{21} & & & \textmd{  and  } & &  & 1
= (f + f^*) 
\end{matrix}$$

\noindent
\underline{PART 2}: Next we show that the CYBE for $r$ (that
is, Equation \ref{YBequation}), translates to a nice expression
in terms of the associated function $f$. This part of the proof
requires some modifications to the proof of the Lie algebra 
case. 

We have:
$$ r = (f \otimes 1) \Omega = \sum_{\alpha}(-1)^{|I_{\alpha}|}
f(I_{\alpha}) \otimes {I_{\alpha}}^* $$

Let us write the three terms of the CYBE:
$$ \begin{array}{llll}
& [r^{12},r^{13}] &=& 
\sum_{\alpha, \beta} (-1)^{|I_{\alpha}| + |I_{\beta}|} 
(-1)^{|I_{\alpha}||I_{\beta}|}
\left [
f(I_{\alpha}) , f(I_{\beta})
\right ]
\otimes {I_{\alpha}}^* \otimes {I_{\beta}}^*  \\ \\
& [r^{12},r^{23}] &=& 
\sum_{\alpha, \beta} (-1)^{|I_{\alpha}| + |I_{\beta}|} 
f(I_{\alpha}) \otimes 
\left [ 
{I_{\alpha}}^*, f(I_{\beta})
\right ]
\otimes {I_{\beta}}^* \\ \\
& [r^{13},r^{23}] 
&=&
\sum_{\alpha, \beta} (-1)^{|I_{\alpha}| + |I_{\beta}|} 
(-1)^{|I_{\alpha}||I_{\beta}|}
f(I_{\alpha}) \otimes f(I_{\beta}) \otimes 
\left [
{I_{\alpha}}^*, {I_{\beta}}^*
\right ] 
\end{array} $$
\noindent
Here we use: 
$$ \begin{array}{cccl}
& [a \otimes b \otimes 1, c \otimes 1 \otimes d] &=&
(-1)^{|b||c|}  [a,c] \otimes b \otimes d \\
\\
& [a \otimes b \otimes 1, 1 \otimes c \otimes d] &=& a \otimes
[b,c] \otimes d \\
\\
& [a \otimes 1 \otimes b, 1 \otimes c \otimes d] &=&
(-1)^{|b||c|} a \otimes c \otimes [b,d] 
\end{array}$$
\noindent
and:
$$ \begin{matrix} |I_{\alpha}| = |{I_{\alpha}}^*| & & &
\textmd{and} & & &  |f(I_{\alpha})| = |I_{\alpha}| \end{matrix}
$$ \noindent 
(We assume $f$ is even.) We rewrite the last sum so that it ends
with $\otimes {I_{\beta}}^*$: 
\begin{eqnarray*}
&& \sum_{\alpha, \beta} (-1)^{|I_{\alpha}| +
|I_{\beta}|}  (-1)^{|I_{\alpha}||I_{\beta}|}
f(I_{\alpha}) \otimes f(I_{\beta}) \otimes 
\left [
{I_{\alpha}}^*, {I_{\beta}}^*
\right ] \\
\\
&=& - \sum_{\alpha , \beta} (-1)^{|I_{\alpha}| + |I_{\beta}|}
f(I_{\alpha}) \otimes f([{I_{\alpha}}^*, I_{\beta}]) 
\otimes {I_{\beta}}^*
\end{eqnarray*}
\noindent
where we use the invariance of the form, and the supersymmetry
of the bracket. 

Therefore we can rewrite the CYBE as: 
$$ \sum_{\alpha,\beta} (-1)^{|I_{\beta}|} 
\begin{pmatrix}
(-1)^{|I_{\alpha}|} (-1)^{|I_{\alpha}||I_{\beta}|}
[f(I_{\alpha}),f(I_{\beta})] \otimes {I_{\alpha}}^* \\
\\
+ (-1)^{|I_{\alpha}|}
f(I_{\alpha}) \otimes [{I_{\alpha}}^*, f(I_{\beta})] \\
\\
- (-1)^{|I_{\alpha}|}
f(I_{\alpha}) \otimes f([{I_{\alpha}}^*, I_{\beta}]) 
\end{pmatrix}
\otimes {I_{\beta}}^* = 0$$

Since the $\{ {I_{\beta}}^*\}$ form a basis for $\g$, this
last equation implies that, for any choice of $\beta$, we have:
$$ \begin{pmatrix}
\sum_{\alpha} (-1)^{|I_{\alpha}|}
(-1)^{|I_{\alpha}||I_{\beta}|}
[f(I_{\alpha}),f(I_{\beta})] \otimes {I_{\alpha}}^* 
\\
\\
+  \sum_{\alpha} (-1)^{|I_{\alpha}|} f(I_{\alpha}) \otimes
[{I_{\alpha}}^*, f(I_{\beta})]  
\\
\\
-  \sum_{\alpha} (-1)^{|I_{\alpha}|}
f(I_{\alpha}) \otimes
f([{I_{\alpha}}^*, I_{\beta}])    
\end{pmatrix} = 0$$

We want to rewrite the second and the third sums so that they
end with $\otimes {I_{\alpha}}^*$. After some calculation, the
second term becomes:
$$ \sum_{\alpha} (-1)^{|I_{\alpha}|} f(I_{\alpha}) \otimes
[{I_{\alpha}}^*, f(I_{\beta})] = \boxed{ - \sum_{\alpha}
(-1)^{|I_{\alpha}|} (-1)^{|I_{\alpha}||I_{\beta}|}
f([I_{\alpha}, f(I_{\beta})]) \otimes {I_{\alpha}}^*} $$
\noindent
The third sum splits into two different sums when we use
Equation \ref{unitarity}: 
\begin{eqnarray*}
&-& \sum_{\alpha} (-1)^{|I_{\alpha}|}
f(I_{\alpha}) \otimes
f([{I_{\alpha}}^*, I_{\beta}]) \\
\\
&=& - \sum_{\alpha} (-1)^{|I_{\alpha}|}
f(I_{\alpha}) \otimes
[{I_{\alpha}}^*, I_{\beta}]
+ \sum_{\alpha} (-1)^{|I_{\alpha}|}
f (f(I_{\alpha})) \otimes
[{I_{\alpha}}^*, I_{\beta}]  
\end{eqnarray*}

\noindent
We calculate these terms separately:
$$ - \sum_{\alpha} (-1)^{|I_{\alpha}|}
f(I_{\alpha}) \otimes [{I_{\alpha}}^*, I_{\beta}] = \boxed{
\sum_{\alpha} (-1)^{|I_{\alpha}|} (-1)^{|I_{\alpha}||I_{\beta}|}
f([I_{\alpha},I_{\beta}]) \otimes {I_{\alpha}}^*} $$
\noindent
and similar but more difficult calculations yield: 
$$ \sum_{\alpha} (-1)^{|I_{\alpha}|}
f (f(I_{\alpha})) \otimes
[{I_{\alpha}}^*, I_{\beta}]  =  \boxed{ - \sum_{\alpha}
(-1)^{|I_{\alpha}|} (-1)^{|I_{\alpha}||I_{\beta}|}
f([f(I_{\alpha}),I_{\beta}] )\otimes {I_{\alpha}}^* } $$
\noindent
Hence we get: 
$$ \sum_{\alpha} (-1)^{|I_{\alpha}|}
(-1)^{|I_{\alpha}||I_{\beta}|}  \begin{pmatrix}
[f(I_{\alpha}),f(I_{\beta})] 
&-& f([I_{\alpha}, f(I_{\beta})]) \\
\\
+ f([I_{\alpha},I_{\beta}])
&-& f([f(I_{\alpha}),I_{\beta}] )
\end{pmatrix} 
\otimes {I_{\alpha}}^* = 0 $$
\noindent
Again using the fact that the $\{ {I_{\alpha}}^*\}$ form a basis
for $\g$, we obtain, for all $\alpha,\beta$:
\begin{equation*}
[f(I_{\alpha}),f(I_{\beta})] - f([I_{\alpha},f(I_{\beta})]) +
f([I_{\alpha},I_{\beta}])  - f([f(I_{\alpha}),I_{\beta}]) = 0
\end{equation*} 
\begin{eqnarray*} 
&\Rightarrow& [f(I_{\alpha}),f(I_{\beta})] =
f([I_{\alpha},f(I_{\beta})]) - f([I_{\alpha},I_{\beta}])  +
f([f(I_{\alpha}),I_{\beta}])  \\ 
\\
&\Rightarrow& (f-1)[f(I_{\alpha}),f(I_{\beta})] =
f([(f-1)(I_{\alpha}),(f-1)(I_{\beta})]), 
\end{eqnarray*}

\noindent
which is equivalent to
\begin{equation*}
\tag{\ref{YBforf}}
\boxed{
(f-1)[f(x),f(y)] = f([(f-1)(x),(f-1)(y)]) \textmd{   for all }
x, y \in \g}
\end{equation*}

\noindent
This proves one direction of the lemma. To see the other
direction, we need only trace the steps above  backwards. Hence
one can easily see that a function $f$ satisfying  Equations
\ref{unitarityforf} and \ref{YBforf} will correspond to an 
$r-$matrix ${r \in \g \otimes \g}$ that satisfies Equations 
\ref{unitarity} and \ref{YBequation}.  This completes the
proof. $\blacksquare$
\\
\\
Writing $r_0 \in \h \otimes \h$ as $r_0 = (f_0 \otimes 1)
\Omega_0$ for a linear map $f_0 : \h \rightarrow \h$, we have:

\begin{lemma}
The system:
\begin{equation*}
\tag{\ref{unitarityinh}}
{r_0}^{12} + {r_0}^{21} = \Omega_0 
\end{equation*}
\begin{equation*}
\tag{\ref{YBinh}}
(\tau(\alpha) \otimes 1) (r_0) + (1 \otimes \alpha) (r_0) = 0
\textmd{  for  all } \alpha \in \Gamma_1 
\end{equation*}
\noindent
is equivalent to the system:
\begin{equation}
\label{unitarityforf0}
f_0 + {f_0}^* = 1
\end{equation}
\begin{equation}
\label{YBforf0}
f_0(h_{\alpha}) = (f_0- 1) (h_{\tau(\alpha)}) 
\textmd{  for  all } \alpha \in \Gamma_1 
\end{equation}
\end{lemma}

\noindent
\textbf{Proof:} We will prove a stronger result. Namely, we will
prove that, for any ${1 \le s,t \le r}$, the equations:
$$\begin{matrix} 
{r_0}^{12} + {r_0}^{21} = \Omega_0 \\
\\
(\alpha_t \otimes 1) (r_0) + (1 \otimes \alpha_s) (r_0) = 0
\end{matrix}$$
\noindent
are equivalent to the equations:
$$\begin{matrix}
f_0 + {f_0}^* = 1 \\ \\
f_0(h_{\alpha_s}) = (f_0- 1) (h_{\alpha_t}) 
\end{matrix} $$

\noindent
It is easy to see the equivalence of the first equations:
\begin{eqnarray*}
{r_0}^{12} + {r_0}^{21} &=& (f_0 \otimes 1 + 1 \otimes f_0)
\Omega_0 \\
&=& ((f_0 + {f_0}^*) \otimes 1) \Omega_0 \\
&=& \Omega_0 \textmd{  if and only if } f_0 + {f_0}^* = 1
\end{eqnarray*}

\noindent
Next we look at $(\alpha_t \otimes 1) r_0 + (1 \otimes
\alpha_s) r_0$. This is equal to:
\begin{eqnarray*}
&& (\alpha_t \otimes 1) 
\left ( 
\sum_i f_0(h_{\alpha_i}) \otimes {h_{\alpha_i}}^* 
\right ) + (1 \otimes \alpha_s) 
\left ( 
\sum_i f_0(h_{\alpha_i}) \otimes {h_{\alpha_i}}^* 
\right ) \\
\\
&=& \sum_i \alpha_t(f_0(h_{\alpha_i})) \cdot {h_{\alpha_i}}^* 
+ \sum_i \alpha_s({f_0}^*({h_{\alpha_i}}^*)) \cdot
h_{\alpha_i} 
\\
\\
&=& \sum_i \alpha_t 
\left (
\sum_k (f_0(h_{\alpha_i}),h_{\alpha_k}) {h_{\alpha_k}}^* 
\right ) \cdot {h_{\alpha_i}}^* 
+ \sum_i \alpha_s
\left (
\sum_k ({f_0}^*({h_{\alpha_i}}^*),h_{\alpha_k}) {h_{\alpha_k}}^* 
\right ) \cdot h_{\alpha_i} \\
\\
&=& \sum_{i,k} (f_0(h_{\alpha_i}),h_{\alpha_k})
\alpha_t({h_{\alpha_k}}^*) \cdot {h_{\alpha_i}}^* 
+ \sum_{i,k}  ({f_0}^*({h_{\alpha_i}}^*),h_{\alpha_k}) 
\alpha_s({h_{\alpha_k}}^*) \cdot h_{\alpha_i} 
\end{eqnarray*}

\noindent
We have: 
$$ \alpha_s({h_{\alpha_k}}^*) = (h_{\alpha_s}, {h_{\alpha_k}}^*)
= \delta_{s,k}  \; \; \; \; \; \; and \; \; \; \; \; \;
\alpha_t({h_{\alpha_k}}^*) = (h_{\alpha_t}, {h_{\alpha_k}}^*) =
\delta_{t,k}$$ \noindent
Therefore the above expression becomes:
\begin{eqnarray*} 
&& \sum_i (f_0(h_{\alpha_i}),h_{\alpha_t}) {h_{\alpha_i}}^*
+ \sum_i ({f_0}^*({h_{\alpha_i}}^*), h_{\alpha_s})
h_{\alpha_i} \\
\\
&=& \sum_i (h_{\alpha_i}, {f_0}^*(h_{\alpha_t})) {h_{\alpha_i}}^*
+ \sum_i ({h_{\alpha_i}}^*, f_0(h_{\alpha_s}))
h_{\alpha_i} \\
\\
&=& {f_0}^*(h_{\alpha_t}) + f_0(h_{\alpha_s}) \\
\\
&=& (1- f_0)(h_{\alpha_t}) + f_0(h_{\alpha_s})
\end{eqnarray*}

\noindent
Summarizing, we have shown that:
$$ (\alpha_t \otimes 1) r_0 + (1 \otimes \alpha_s) r_0 = 
(1- f_0)(h_{\alpha_t}) + f_0(h_{\alpha_s}) $$
\noindent
Hence one is equal to zero if and only if the other is. This completes the
proof of the lemma. $\blacksquare$
\\
\\
We also need:

\begin{lemma} 
The system of equations
\begin{equation*}
\tag{\ref{unitarityinh}}
{r_0}^{12} + {r_0}^{21} = \Omega_0 
\end{equation*}
\begin{equation*}
\tag{\ref{YBinh}}
(\tau(\alpha) \otimes 1) (r_0) + (1 \otimes \alpha) (r_0) = 0
\textmd{  for  all } \alpha \in \Gamma_1 
\end{equation*}
\noindent
is consistent.
\end{lemma}

\begin{remark}
The arguments used to prove this lemma are the same as for
the Lie algebra case (see \cite{BD2} for details), and hence will
not be included here.  
\end{remark}

These lemmas allow us to translate the conditions (Equations 
\ref{unitarityinh} and \ref{YBinh}) on the continuous
parameter used in  the main theorem into conditions on a linear
map ${f_0 : \h \rightarrow  \h}$. Also recall that Lemma
\ref{rstatementsintofstatements} translates the CYBE and the
unitarity condition (i.e. Equations \ref{unitarity} and 
\ref{YBequation}) into conditions on the associated linear map
$f : \g  \rightarrow \g$, namely Equations \ref{unitarityforf}
and \ref{YBforf}.

So from now on, we will be using linear maps $f, f_0$
and so on, interchangeably with their $2-$tensor equivalents 
$r, r_0$ etc. Then, we can reformulate our problem in
the following manner: Given an admissible triple
$(\Gamma_1,\Gamma_2,\tau)$ with a linear map $f_0 : \h
\rightarrow \h$ satisfying Equations \ref{unitarityforf0} and
\ref{YBforf0}, construct a linear map  ${f : \g
\rightarrow \g}$ satisfying Equations \ref{unitarityforf} and 
\ref{YBforf}.

\section{The Cayley Transform}

We will now consider a variation on the theme of Cayley
transforms. For the time being assume that we have a linear
function $f : \g \rightarrow \g$ with $(f-1)$ invertible. Then
\textbf{the Cayley transform of} $f$ is $\Theta =
\frac{f}{f-1}$. Then the adjoint of this function will be
$\Theta^* = \frac{f^*}{(f-1)^*} = \frac{1-f}{-f}$ if our $f$
satisfies Equation \ref{unitarityforf}: 
\begin{equation*} 
\tag{\ref{unitarityforf}} 
f + f^* = 1
\end{equation*}

Then one can see that $\Theta \Theta^* = 1$, and so $\Theta$
preserves the invariant form. If
we also assume that $f$ satisfies Equation \ref{YBforf}:
\begin{equation*} \tag{\ref{YBforf}}
(f-1) ([f(x),f(y)]) = f([(f-1)(x),(f-1)(y)])
\end{equation*}
\noindent
we will get:
\begin{equation*}
[\Theta(x),\Theta(y)] = \Theta([x,y])
\end{equation*}
\noindent
which implies that $\Theta$ is a Lie superalgebra automorphism.

However for the Lie superalgebras we care about, i.e. for simple Lie
superalgebras, the above will not work. To see this, assume that $f$
is a linear map satisfying Equations \ref{unitarityforf} and 
\ref{YBforf}, ${f-1}$ is invertible, and $\Theta$ is defined as above. 
Then look at ${\Theta - 1}$. This is given by ${\frac{f}{f-1} - 1 }= 
\frac{1}{f-1}$, i.e. it is the inverse of ${f-1}$. So we must have that 
${det(\Theta - 1) \ne 0}$. But this is not allowed for a simple Lie 
superalgebra:

\begin{lemma}
\label{Cayleylemma1}
If $\Theta$ is an automorphism of a finite dimensional (classical)
simple Lie superalgebra, then ${det(\Theta - 1) = 0}$.
\end{lemma}

\noindent
\textbf{Proof:} We will need the following result from \cite{BD2} (this 
is Theorem 9.2 there):
\\
\\
\textit{If $\phi$ is an automorphism of a semisimple Lie algebra
$L$, then there exists a nonzero element $x \in L$ such that
$\phi(x) = x$.}
\\
\\
Now the automorphism $\Theta$ of a simple Lie superalgebra $\g$ restricts
to a (Lie algebra) automorphism $\theta$ on $\g_{\overline{0}}$,
the even part of $\g$. $\g_{\overline{0}}$ is reductive with
nontrivial ${\g_{\overline{0}}}' = {[\g_{\overline{0}},
\g_{\overline{0}}]}$  \footnotemark. 
\footnotetext{ The even part of a classical simple Lie
superalgebra will be reductive and will decompose into a direct
sum of its derived algebra (which is nonempty) and some abelian
subalgebra.} 
${\g_{\overline{0}}}'$ is semisimple and $\theta$ restricts to
an automorphism $\phi$ on ${\g_{\overline{0}}}'$. But then the
above result gives us some nonzero ${x \in
{\g_{\overline{0}}}'}$ with $\phi(x) = x$.  Going back to our
$\Theta$ we see that $\Theta(x) = x$ and hence ${x \in
Ker(\Theta - 1)}$. Thus we must have that ${det(\Theta - 1) =
0}$. This proves the lemma. $\blacksquare$  
\\  
\\ 
\noindent
Thus Equations \ref{unitarityforf} and \ref{YBforf} will imply
that $f-1$  is not invertible, as the invertibility assumption
leads to a  contradiction with Lemma
\ref{Cayleylemma1} \footnotemark.  Therefore, we cannot define
the Cayley transform as above for the functions we are
interested in. 
\footnotetext{Using Equation \ref{unitarityforf} once again, we
see that  $f$ is not invertible, either. Thus any solution $r$
to the system of  Equations \ref{unitarityforf} and \ref{YBforf}
will be degenerate.}

However it turns out that we can modify our definition and still get a
lot of what we want: First note that for any linear operator $f$,
${Ker(f) \subset Im(f-1)}$ and ${Ker(f-1) \subset Im(f)}$. Then we define
the \textbf{Cayley transform of} $f$ to be the function {$\Theta :
Im(f-1)/Ker(f) \rightarrow Im(f)/Ker(f-1)$} that maps
$\overline{(f-1)(x)}$ to $\overline{f(x)}$. [It is easy to check that
this is well-defined.] This version of the Cayley transform will be
sufficient for our purposes. We have:

\begin{lemma}
\label{Cayleylemma2}
Let $f : \g \rightarrow \g$ be a linear map satisfying:
\begin{equation*}
\tag{\ref{unitarityforf}}
f + f^* = 1
\end{equation*}
\noindent
Then $Ker(f) = Im(f-1)^{\perp}$, $Ker(f-1) = Im(f)^{\perp}$, and the map
$\Theta$ preserves the invariant form. Furthermore, $f$ satisfies:
\begin{equation*}
\tag{\ref{YBforf}}
(f-1)[f(x),f(y)] = f([(f-1)(x),(f-1)(y)]),
\end{equation*}
\noindent
if and only if $Im(f)$ and $Im(f-1)$ are Lie subsuperalgebras of $\g$, and
$\Theta$ is a Lie superalgebra isomorphism.
\end{lemma}

\begin{remark}
This lemma will be valid for any simple Lie superalgebra $\g$ with a
non-degenerate Killing form. In this case its proof will be
exactly the same as the proof of the analogous result in the
Lie algebra case. See \cite{BD2}. 
\end{remark}

\section{The Construction - End of the Proof of the Theorem}

For a given admissible triple $(\Gamma_1,\Gamma_2, \tau)$, and
a linear map  ${f_0  : \h \rightarrow \h}$ satisfying Equations
\ref{unitarityforf0} and \ref{YBforf0}, we want to construct a function 
${f: \g \rightarrow \g}$ that will satisfy Equations \ref{unitarityforf} 
and \ref{YBforf}. Here is how we proceed:

For the admissible triple $(\Gamma_1, \Gamma_2, \tau)$, we
define $\overline{\Gamma}_i$ and $\overline{\tau}$ as before.
Also we define the following Lie subsuperalgebras of $\g$:   
$$ \begin{matrix}  
\h_i &=& \bigoplus_{\alpha \in \Gamma_i} \C h_{\alpha} & & & & & 
\g_i = \h_i \oplus \sum_{\alpha \in \overline{\Gamma}_i}   
\left (
\g_{\alpha} \oplus \g_{-\alpha}
\right ) 
\end{matrix}$$
\noindent
and:
$$ \begin{matrix}
{\mathfrak{n}_i}^+ &=& \sum_{\alpha \in \Delta^+ /
\overline{\Gamma}_i} \g_{\alpha} & & & & & 
{\mathfrak{p}_i}^+ &=& \g_i + {\mathfrak{n}_i}^+ \\
\\
{\mathfrak{n}_i}^- &=& \sum_{\alpha \in \Delta^+ /
\overline{\Gamma}_i} \g_{-\alpha} & & & & & 
{\mathfrak{p}_i}^- &=& \g_i + {\mathfrak{n}_i}^- 
\end{matrix}
$$
\noindent
We can see that the ${\mathfrak{n}_i}^{+/-}$ are ideals
in ${\mathfrak{p}_i}^{+/-}$. 
\\
\\
Next if $f_0 : \h \rightarrow \h$ satisfies:
\begin{equation*}
\tag{\ref{YBforf0}}
f_0(h_{\alpha}) = (f_0- 1) (h_{\tau(\alpha)}) 
\textmd{  for  all } \alpha \in \Gamma_1 
\end{equation*}
\noindent
we get:
$$ \begin{matrix} h_{\alpha} &=& (f_0 -1)(h_{\tau(\alpha)} -
h_{\alpha}) & & & & & h_{\tau(\alpha)} &=& f_0( h_{\tau(\alpha)} -
h_{\alpha}) 
\end{matrix} $$
\noindent
for all $\alpha \in \Gamma_1$. This implies that $h_{\alpha}
\in Im(f_0 - 1)$ and $h_{\tau(\alpha)} \in Im(f_0)$. Therefore we
have that: $\h_1 \subset Im(f_0 -1)$, and $\h_2 \subset
Im(f_0)$. 
\\
\\
Now fix a Weyl-Chevalley basis $\{ X_{\alpha_i}, Y_{\alpha_i},
H_{\alpha_i} | \alpha_i \in \Gamma\}$. It is known that such a
set of generators exists and satisfies the usual Serre-type
relations, (see \cite{LS} for details).We define a map $\phi$ by 
$$ \phi(X_{\alpha}) = X_{\tau(\alpha)} 
\; \; \; \; \; \; \; 
\phi(Y_{\alpha}) = Y_{\tau(\alpha)}  
\; \; \; \; \; \; \; 
\phi(H_{\alpha}) = H_{\tau(\alpha)} $$ 
\noindent 
for all $\alpha \in \Gamma_1$. Then this can be extended to an
isomorphism ${\phi : \g_1 \rightarrow \g_2}$ because the
relations between $X_{\alpha}, Y_{\alpha}, H_{\alpha}$ for
$\alpha \in \Gamma_1$ will be the same as the relations between 
$X_{\tau(\alpha)}, Y_{\tau(\alpha)}, H_{\tau(\alpha)}$ for
$\alpha \in \Gamma_1$ [Here we are using the fact that $\tau$
is an isometry preserving grading]. Note that $\phi^{-1}$ is a map from
$\g_2$ onto $\g_1$. Since $\tau$ is an isometry, ${(\phi(x),
y)_{\g_2} = (x, \phi^{-1}(y))_{\g_1}}$ for all ${x \in \g_1}, {y
\in \g_2}$. But $\phi^*$ should map $\g_2$ into $\g_1$ and
satisfy exactly the same conditions; hence $\phi^* = \phi^{-1}$.

Now in each root space $\g_{\alpha}$, we can choose an element
$e_{\alpha}$ such that $(e_{\alpha},e_{-\alpha}) =1$ for any
$\alpha$ and $\phi(e_{\alpha}) = e_{\overline{\tau}(\alpha)}$
for all $\alpha \in \overline{\Gamma}_1$. The fact that there
are no cycles for $\tau$ will ensure that such a choice is
possible. 

Next define a linear map as follows: 
$$ \psi(x) = 
\left \{ 
\begin{matrix} 
\phi(x) && \textmd{if } x \in \g_1 \\ 
0 && \textmd{if }  x \in {\mathfrak{n}_1}^+
\end{matrix}
\right . $$ 
\noindent
This restricts to a map on 
$\mathfrak{n}_+ = \bigoplus_{\alpha > 0} \g_{\alpha}$, since
${\mathfrak{n}_+ = (\g_1 \cap \mathfrak{n}_+) \oplus 
{\mathfrak{n}_1}^+}$. The proof of the following lemma is
exactly the same as in the Lie algebra case (see \cite{BD2}):

\begin{lemma} 
det$(\psi - 1)$ is nonzero if and only if $\tau$ satisfies the
second condition in the definition of an admissible triple.
\end{lemma}

Recall that we started with an admissible triple. The above lemma
then says that $\psi -1$ is invertible. Therefore we can define a
function on $\mathfrak{n}_+$ by: 
$$f_+ = \frac{\psi}{\psi - 1} = -(\psi + \psi^2 + \cdots) $$
\noindent  
Clearly the sum on the right hand side is finite as $\psi$ is
nilpotent. Notice that $\psi^*$ and so ${f_+}^*$ are maps on
$\mathfrak{n}_-= \bigoplus_{\alpha < 0} \g_{\alpha}$, since the
Killing form induces a non-degenerate pairing of $\mathfrak{n}_+$
with $\mathfrak{n}_-$. 

Now define a linear map on $\mathfrak{n}_-$ by:
$$ f_- = 1 - {f_+}^* = 1 + \psi^* + {\psi^*}^2 + \cdots $$
\noindent  
Then define $f$ to be the function whose restriction to
$\h$, $\mathfrak{n}_+$, $\mathfrak{n}_-$ is $f_0$, $f_+$, $f_-$,
respectively. We have:
\begin{eqnarray*}
f + f^* &=& (f_0 + f_+ + f_-) + (f_0 + f_+ + f_-)^* \\
&=&  (f_0 + {f_0}^*) + (f_+ + {f_-}^*) + ({f_+}^* + f_-) \\
&=&  1_{\h} + 1_{\mathfrak{n}_+} + 1_{\mathfrak{n}_-} \\
&=& 1_{\g}
\end{eqnarray*}

To see that $f$ satisfies Equation \ref{YBforf}, we will use
Lemma \ref{Cayleylemma2}. Recall that for a linear map ${f: \g
\rightarrow \g}$,  we defined the Cayley transform to be the
function 
$${\Theta :  Im(f-1)/Ker(f) \rightarrow Im(f)/Ker(f-1)}$$ 
\noindent
that maps  $\overline{(f-1)(x)}$ to $\overline{f(x)}$. Then we
have seen before that if $f$ satisfies Equation
\ref{unitarityforf} i.e.  ${f + f^* = 1}$, then ${Ker(f) =
Im(f-1)^{\perp}}$, ${Ker(f-1) =  Im(f)^{\perp}}$, and ${\Theta
{\Theta}^* = 1}$. Furthermore, $f$ satisfies Equation
\ref{YBforf} if and only if $Im(f)$ and $Im(f-1)$ are Lie
subsuperalgebras of  $\g$, and $\Theta$ is a Lie superalgebra
isomorphism. 

Thus our problem now reduces to showing that ${C_1 = Im(f-1)}$ and
${C_2 = Im(f)}$ are Lie subsuperalgebras of $\g$, and that the Cayley
transform $\Theta$ of $f$ is a Lie superalgebra isomorphism.

We have 
\begin{eqnarray*}
C_1 &=& Im(f-1) = {Im(f_0-1) \cup Im(f_+-1) \cup Im(f_--1)} \\
C_2 &=& Im(f) = {Im(f_0) \cup Im(f_+) \cup Im(f_-)}
\end{eqnarray*}

We have seen that ${Im(f_0-1) \supset \h_1}$ and ${Im(f_0)
\supset  \h_2}$. We will therefore define $V_1$, $V_2$ as
(vector) subspaces of $\h$ such that ${Im(f_0 -1) = \h_1 \oplus
V_1}$ and ${Im(f_0) = \h_2  \oplus V_2}$. 
\\
\\
\indent
In the Lie algebra case, the Killing form restricts
to a positive definite non-degenerate form on (the real subspace
generated by ${ \{H_{\alpha} | \alpha \in \Gamma \}}$ of) $\h$.
So we can define the orthogonal complements of $\h_1$ and $\h_2$
with respect to this form; call these ${\h_1}^{c}$ and
${\h_2}^{c}$; then we have: $\h = {\h_1 \oplus {\h_1}^c} = {\h_2
\oplus {\h_2}^c}$. Then for a fixed $f_0$, the two subspaces
$V_1$ and $V_2$ are uniquely determined if we add the condition
that ${V_i \subset {\h_i}^{c}}$. 

In the super case, this is no longer possible; the Cartan
subalgebra $\h$ may have isotropic elements and subspaces of
$\h$ may intersect their orthogonal complements nontrivially.
However in our case we still can define ${\h_i}^c$ as follows:
$$ {\h_i}^c = \bigoplus_{\alpha \in \Gamma \backslash \Gamma_i}
\C h_{\alpha} $$
\noindent 
Thus we still can write ${\h = \h_i \oplus {\h_i}^c}$, and still
can demand that ${V_i \subset {\h_i}^c}$. This choice of the
$V_i$ is then again well-defined, but clearly depends on our
choice for $\Gamma$. 

Next we compute:
$$\begin{array}{llll}
&Im(f_+-1) &=& Im(\frac{1}{\psi -1}) = \mathfrak{n}_+ \\
\\
&Im(f_--1) &=& Im(\frac{\psi^*}{1 - \psi^*}) = Im(\psi^*) = \g_1
\cap \mathfrak{n}_- \\
\\ 
&Im(f_+) &=& Im(\frac{\psi}{1-\psi}) = Im(\psi) = \g_2 \cap
\mathfrak{n}_+ \\
\\
&Im(f_-) &=&  Im(\frac{1}{1 - \psi^*}) = \mathfrak{n}_-
\end{array}$$
\noindent
where we use the fact that $\psi -1$ is invertible. The above
then yields: 
$$ \begin{matrix}  
C_1 = {\mathfrak{p}_1}^+ \oplus V_1 & & &
C_2 = {\mathfrak{p}_2}^- \oplus V_2 
\end{matrix} $$
\noindent
It is now easy to check that $C_1$ and $C_2$ are both closed under the 
bracket and hence are Lie subsuperalgebras of $\g$.

Finally we need to see that the Cayley transform $\Theta$ is a Lie 
superalgebra isomorphism. Now we note that by the last lemma above, 
${C_i \supset {C_i}^{\perp}}$. So we have:
$$ {C_1}^{\perp} = ({\mathfrak{p}_1}^+ \oplus V_1)^{\perp} = 
{\mathfrak{n}_1}^+ \oplus (\h_1 \oplus V_1)^{\perp} =
{\mathfrak{n}_1}^+ \oplus ({\h_1}^{\perp} \cap {V_1}^{\perp}) \subset 
{\mathfrak{p}_1}^+ 
\oplus V_1$$
\noindent
and similarly:
$$ {C_2}^{\perp} = ({\mathfrak{p}_2}^- \oplus V_2)^{\perp} = 
{\mathfrak{n}_2}^- \oplus (\h_2 \oplus V_2)^{\perp}  = 
{\mathfrak{n}_2}^- \oplus ({\h_2}^{\perp} \cap {V_2}^{\perp}) \subset 
{\mathfrak{p}_2}^- \oplus V_2$$
\noindent
Hence we have: 
$$ {\h_i}^{\perp} \cap {V_i}^{\perp} \subset \h_i \oplus V_i $$
\noindent
which gives us: 
$$ C_1 / {C_1}^{\perp} = \frac{{\mathfrak{p}_1}^+ \oplus V_1 }  
{({\mathfrak{p}_1}^+ \oplus V_1)^{\perp}} = 
\left (
\bigoplus_{\alpha \in 
\overline{\Gamma}_1} \g_{\alpha} \oplus \g_{-\alpha}
\right ) 
\oplus \frac {\h_1 \oplus V_1}  {{\h_1}^{\perp} \cap {V_1}^{\perp}} $$

\noindent
and similarly:
$$ C_2 / {C_2}^{\perp} = \frac{{\mathfrak{p}_2}^- \oplus V_2 }  
{({\mathfrak{p}_2}^- \oplus V_2)^{\perp}} = 
\left (
\bigoplus_{\alpha \in 
\overline{\Gamma}_2} \g_{\alpha} \oplus \g_{-\alpha}
\right ) 
\oplus \frac {\h_2 \oplus V_2}  {{\h_2}^{\perp} \cap {V_2}^{\perp}} $$
\noindent
We have already seen that the $C_i$ are Lie subsuperalgebras.
Since  ${C_i}^{\perp}$ is an ideal, we have a Lie superalgebra
structure on  $C_i / {C_i}^{\perp}$. But ${[\g_{\alpha},
\g_{-\alpha}] = \C  H_{\alpha}}$, therefore we must have a
complete copy of $\h_i$ and so  a copy of $\g_i$ in $C_i /
{C_i}^{\perp}$. This implies that  ${{\h_i}^{\perp} \cap
{V_i}^{\perp}}$ lies in $V_i$, and we have:  $$ C_i /
{C_i}^{\perp} = \g_i \oplus \frac{V_i} {{\h_i}^{\perp} \cap 
{V_i}^{\perp}} $$

\noindent
So we need to show that:
$$\Theta : \g_1 \oplus \frac {V_1}  {{\h_1}^{\perp} \cap {V_1}^{\perp}} 
\rightarrow \g_2 \oplus \frac {V_2}  {{\h_2}^{\perp} \cap {V_2}^{\perp}} 
$$
\noindent
is a Lie superalgebra isomorphism. 

We first note that $\Theta(x) = \phi(x)$ for all $x \in 
\g_1$. Indeed if $\alpha \in \Gamma_1$, we have:
$$ X_{\alpha} = (f_+ -1) (X_{\tau(\alpha)} - X_{\alpha}) $$ 
\noindent
and so is mapped via $\Theta$ to:
$$ f_+ (X_{\tau(\alpha)} - X_{\alpha}) = X_{\tau(\alpha)} $$
\noindent
And similarly:
$$ Y_{\alpha} = (f_- -1) (Y_{\tau(\alpha)} - Y_{\alpha}) $$ 
\noindent
is mapped via $\Theta$ to:
$$ f_- (Y_{\tau(\alpha)} - Y_{\alpha}) = Y_{\tau(\alpha)} $$
\noindent
Also it is easy to see that since ${H_{\alpha} = (f_0 - 1) 
(H_{\tau(\alpha)} - H_{\alpha})}$ for each $\alpha \in \Gamma_1$, 
$\Theta$ sends $H_{\alpha}$ to ${f_0(H_{\tau(\alpha)} -
H_{\alpha})} =  H_{\tau(\alpha)}$. Hence the restriction of
$\Theta$ to $\g_1$ is exactly the Lie superalgebra isomorphism
$\phi$. 

Next we look at what $\Theta$ does on the Cartan part of the 
$C_i/{C_i}^{\perp}$. We have that: 
$$ C_i/{C_i}^{\perp} = 
\left (
\bigoplus_{\alpha \in
\overline{\Gamma}_i} \g_{\alpha} \oplus \g_{-\alpha}
\right )
\oplus \frac {\h_i \oplus V_i}  {{\h_i}^{\perp} \cap {V_i}^{\perp}} $$
\noindent
i.e. we can rewrite $C_i/{C_i}^{\perp}$ as the direct sum of a
Cartan part and a non-Cartan part. Then the previous arguments
show that $\Theta$ maps the non-Cartan part of ${C_1 /
{C_1}^{\perp}}$ into the non-Cartan part of ${C_2 /
{C_2}^{\perp}}$ like $\phi$ does. And then since $\Theta$
preserves the invariant form, it maps the Cartan part of
$C_1/{C_1}^{\perp}$ to the Cartan part of $C_2/{C_2}^{\perp}$:  
$$ \Theta  \left (({\textmd{non-Cartan of
}C_1/{C_1}^{\perp}})^{\perp}  \right ) = ({\textmd{non-Cartan of
}C_2/{C_2}^{\perp}})^{\perp} $$  
\noindent 
In other words: 
$$
\Theta  \left ( \frac {\h_1 \oplus V_1}  {{\h_1}^{\perp} \cap 
{V_1}^{\perp}} \right ) =  \left ( \frac {\h_2 \oplus V_2} 
{{\h_2}^{\perp} \cap  {V_2}^{\perp}} \right ) $$
\noindent
Since $\frac {\h_i \oplus V_i} {{\h_i}^{\perp} \cap
{V_i}^{\perp}}$ are abelian, $\Theta$ restricts to an
isomorphism there as well. Therefore  $\Theta$ is an
isomorphism. This then will imply that the associated linear
map $f$ satisfies Equations \ref{unitarityforf} and \ref{YBforf}
and so corresponds to an $r-$ matrix satisfying Equations
\ref{unitarity} and \ref{YBequation}.

At this point one needs to check if the function $f$
constructed in this way will yield the tensor $r$ of Equation
$(*)$. This is reasonably straightforward. Hence we have proved
our theorem. $\blacksquare$

\section{Examples: $r-$matrices on $sl(2,1)$}
\label{Sectiononsl21}
 
Recall that two Dynkin diagrams of a given Lie 
superalgebra are not necessarily isomorphic, but one can be
obtained from the other via a chain of odd reflections. 
Therefore while listing all possible admissible triples for a
given Lie superalgebra, we need to take into consideration
all possible Dynkin diagrams. This raises a new
question as to how r-matrices obtained from two nonisomorphic
Dynkin diagrams are related, if at all. We will see that at
least in the case of $sl(2,1)$, if $r$ is the standard r-matrix
associated to a fixed Dynkin diagram $D$, and $D'$ is the Dynkin
diagram obtained from $D$ by the odd reflection
$\sigma_{\alpha}$, then $r'$, the standard r-matrix associated
to $D'$, will be the image of $r$ under the same reflection
$\sigma_{\alpha}$. 

\subsection{Dynkin Diagrams of sl(2,1)}

The roots of $sl(2,1)$ are 
$$ \Delta_{\overline{0}} = \{ \epsilon_1 - \epsilon_2,
\epsilon_2 - \epsilon_1 \}
\; \; \; \; \; \; 
\Delta_{\overline{1}} = \{ \epsilon_1 - \lambda_1, 
\epsilon_2 - \lambda_1, \lambda_1 - \epsilon_1, 
\lambda_1 - \epsilon_2  \} $$
\noindent
where $\epsilon_i$ is the (restriction to the Cartan subalgebra
of $sl(2,1)$ of the) standard basis: ${\epsilon_i(E_{jk}) =
\delta_{i,j} \delta_{i,k}}$, and  ${\lambda_1 = \epsilon_3}$. We
will denote the set of simple roots by $\Gamma$. 

There are six possible Dynkin diagrams:
\begin{enumerate}
\item $\Gamma(D_1) = \{ \epsilon_1 - \epsilon_2, 
\epsilon_2 - \lambda_1 \} $. We will set 
${\alpha_1 = \epsilon_1 - \epsilon_2}$ and ${\alpha_2 = 
\epsilon_2 - \lambda_1}$.  $\alpha_1$ is even, while $\alpha_2$
is odd. The third positive root in this case is ${\alpha_1 +
\alpha_2}$ which is odd. 

\item $\Gamma(D_2) = \{ \epsilon_1 - \lambda_1, 
\lambda_1 - \epsilon_2 \}$. Note that these two roots are
actually ${\alpha_1 + \alpha_2}$ and $-\alpha_2$, and it is easy
to see that $D_2$ is obtained from $D_1$ via the odd reflection
$\sigma_{\alpha_2}$ associated to the root $\alpha_2$. We can
obtain $D_1$ back from $D_2$ by $\sigma_{-\alpha_2}$. Note also
that the third positive root in this case will be $\alpha_1$
which is even. 

\item $\Gamma(D_3) = \{ \lambda_1 - \epsilon_1, 
\epsilon_1 - \epsilon_2 \}$ $ = \{ -\alpha_1 -\alpha_2, \alpha_1
\}$. We note that $D_3$ is obtained from $D_2$ via the odd
reflection $\sigma_{\alpha_1 + \alpha_2}$. Applying
$\sigma_{-\alpha_1 - \alpha_2}$ to $D_3$ will return $D_2$ as
expected. Note also that the third positive root will be
$-\alpha_2$ which is odd.

\item $\Gamma(D_4) = -\Gamma(D_1) = \{ -\alpha_1, -\alpha_2 \}$.
The third positive root in this case will be ${-\alpha_1
-\alpha_2}$ which is odd. 

\item $\Gamma(D_5) = -\Gamma(D_2) = \{ -\alpha_1 - \alpha_2,
\alpha_2 \}$. The third positive root in this case will be
$-\alpha_1$ which is even. Note that $D_5$ is obtained from
$D_4$ via $\sigma_{-\alpha_2}$, and that applying
$\sigma_{\alpha_2}$ to $D_5$ will yield $D_4$ as expected. 

\item $\Gamma(D_6) = -\Gamma(D_3) = \{ \alpha_1 + \alpha_2,
-\alpha_1 \}$. The third positive root will be $\alpha_2$ which
is odd. The odd reflections $\sigma_{-\alpha_1
-\alpha_2}$ and $\sigma_{\alpha_1 + \alpha_2}$ will map $D_5$ and
$D_6$ into one another. 

\end{enumerate}

Hence, up to sign, there are three Dynkin diagrams, and these
can be obtained from one another via a chain of odd reflections
(which change the signs of some of the odd roots but a positive
even root stays positive). Also note that except for the two
Dynkin diagrams $D_2$ and $D_5$, the diagrams consist of one
white circle and one gray circle (standing for two roots of
different parities), and so these diagrams will not allow any
nontrivial admissible triples. Therefore the
construction of our theorem will only yield standard r-matrices
for these cases. In $D_2$ and $D_5$, both simple roots
are odd, and we can actually consider a nontrivial admissible
triple in these cases. Therefore the theorem will give us one
standard r-matrix and two nonstandard r-matrices for each of the diagrams 
$D_2$ and $D_5$.

\subsection{The Standard r-matrices}
\label{standardsubsection}

By construction, given $r_0 \in \h \otimes \h$ satisfying ${r_0 +
{r_0}^{21} = \Omega_0}$ \footnotemark, the standard r-matrix for a 
fixed Dynkin diagram is 
$$r = r_0 + \sum_{\alpha > 0} e_{-\alpha} \otimes e_{\alpha}.$$

\footnotetext{$\Omega_0$ is the Cartan part of the quadratic Casimir 
element $\Omega$ of $\g$.} 

\noindent
So fixing $r_0$ we write down the standard r-matrices for the above 
diagrams: 

\begin{enumerate} 

\item  For $D_1$ the positive roots are $\alpha_1$, $\alpha_2$
and ${\alpha_1 + \alpha_2}$. We let:
$$ e_{\alpha_1} = E_{12}, \;\;\;\;\;  e_{\alpha_2} = E_{23},
\;\;\;\;\; e_{\alpha_1 + \alpha_2} = E_{13}$$
\noindent
Then we choose $e_{-\alpha}$ by ${(e_{\alpha},e_{-\alpha}) = 1}$:
$$ e_{-\alpha_1} = E_{21}, \;\;\;\;\;  e_{-\alpha_2} = E_{32},
\;\;\;\;\; e_{-\alpha_1 - \alpha_2} = E_{31}$$
\noindent
Therefore we get:
$$ r_{st}(D_1) = r_0 + (E_{21} \otimes E_{12}) + (E_{32} \otimes
E_{23})  + (E_{31} \otimes E_{13}) $$

\item For $D_2$ the positive roots are $\alpha_1$, $-\alpha_2$
and ${\alpha_1 + \alpha_2}$. We let:
$$ e_{\alpha_1} = E_{12}, \;\;\;\;\;  e_{-\alpha_2} = E_{32},
\;\;\;\;\; e_{\alpha_1 + \alpha_2} = E_{13}$$
\noindent
Then we choose $e_{-\alpha}$ by ${(e_{\alpha},e_{-\alpha}) = 1}$:
$$ e_{-\alpha_1} = E_{21}, \;\;\;\;\;  e_{\alpha_2} = -E_{23},
\;\;\;\;\; e_{-\alpha_1 - \alpha_2} = E_{31}$$
\noindent
Therefore we get:
$$ r_{st}(D_2) = r_0 + (E_{21} \otimes E_{12}) - (E_{23} \otimes
E_{32})  + (E_{31} \otimes E_{13}) $$

\item For $D_3$ the positive roots are $\alpha_1$, $-\alpha_2$
and ${-\alpha_1 - \alpha_2}$. We let:
$$ e_{\alpha_1} = E_{12}, \;\;\;\;\;  e_{-\alpha_2} = E_{32},
\;\;\;\;\; e_{-\alpha_1 - \alpha_2} = E_{31}$$
\noindent
Then we choose $e_{-\alpha}$ by ${(e_{\alpha},e_{-\alpha}) = 1}$:
$$ e_{-\alpha_1} = E_{21}, \;\;\;\;\;  e_{\alpha_2} = -E_{23},
\;\;\;\;\; e_{\alpha_1 + \alpha_2} = -E_{13}$$
\noindent
Therefore we get:
$$ r_{st}(D_3) = r_0 + (E_{21} \otimes E_{12}) - (E_{23} \otimes
E_{32}) - (E_{13} \otimes E_{31}) $$

\item For $D_4$ the positive roots are $-\alpha_1$, $-\alpha_2$
and ${-\alpha_1 - \alpha_2}$. We let:
$$ e_{-\alpha_1} = E_{21}, \;\;\;\;\;  e_{-\alpha_2} = E_{32},
\;\;\;\;\; e_{-\alpha_1 - \alpha_2} = E_{31}$$
\noindent
Then we choose $e_{-\alpha}$ by ${(e_{\alpha},e_{-\alpha}) = 1}$:
$$ e_{\alpha_1} = E_{12}, \;\;\;\;\;  e_{\alpha_2} = -E_{23},
\;\;\;\;\; e_{\alpha_1 + \alpha_2} = -E_{13}$$
\noindent
Therefore we get:
$$ r_{st}(D_4) = r_0 + (E_{12} \otimes E_{21}) - (E_{23} \otimes
E_{32}) - (E_{13} \otimes E_{31}) $$

\item For $D_5$ the positive roots are $-\alpha_1$, $\alpha_2$
and ${-\alpha_1 - \alpha_2}$. We let:
$$ e_{-\alpha_1} = E_{21}, \;\;\;\;\;  e_{\alpha_2} = E_{23},
\;\;\;\;\; e_{-\alpha_1 - \alpha_2} = E_{31}$$
\noindent
Then we choose $e_{-\alpha}$ by ${(e_{\alpha},e_{-\alpha}) = 1}$:
$$ e_{\alpha_1} = E_{12}, \;\;\;\;\;  e_{-\alpha_2} = E_{32},
\;\;\;\;\; e_{\alpha_1 + \alpha_2} = -E_{13}$$
\noindent
Therefore we get:
$$ r_{st}(D_5) = r_0 + (E_{12} \otimes E_{21}) + (E_{32} \otimes
E_{23}) - (E_{13} \otimes E_{31}) $$

\item For $D_6$ the positive roots are $-\alpha_1$, $\alpha_2$
and ${\alpha_1 + \alpha_2}$. We let:
$$ e_{-\alpha_1} = E_{21}, \;\;\;\;\;  e_{\alpha_2} = E_{23},
\;\;\;\;\; e_{\alpha_1 + \alpha_2} = E_{13}$$
\noindent
Then we choose $e_{-\alpha}$ by ${(e_{\alpha},e_{-\alpha}) =
1}$: $$ e_{\alpha_1} = E_{12}, \;\;\;\;\;  e_{-\alpha_2} =
E_{32}, \;\;\;\;\; e_{-\alpha_1 - \alpha_2} = E_{31}$$
\noindent
Therefore we get:
$$ r_{st}(D_6) = r_0 + (E_{12} \otimes E_{21}) + (E_{32} \otimes
E_{23})  + (E_{31} \otimes E_{13}) $$

\end{enumerate}

Summarizing we have:

\begin{enumerate}
\item $ r_{st}(D_1) = r_0 + (E_{21} \otimes E_{12}) + (E_{32}
\otimes E_{23})  + (E_{31} \otimes E_{13}) $
\item $ r_{st}(D_2) = r_0 + (E_{21} \otimes E_{12}) - (E_{23} \otimes
E_{32})  + (E_{31} \otimes E_{13}) $
\item $ r_{st}(D_3) = r_0 + (E_{21} \otimes E_{12}) - (E_{23} \otimes
E_{32}) - (E_{13} \otimes E_{31}) $
\item $ r_{st}(D_4) = r_0 + (E_{12} \otimes E_{21}) - (E_{23} \otimes
E_{32}) - (E_{13} \otimes E_{31}) $
\item $ r_{st}(D_5) = r_0 + (E_{12} \otimes E_{21}) + (E_{32} \otimes
E_{23}) - (E_{13} \otimes E_{31}) $
\item $ r_{st}(D_6) = r_0 + (E_{12} \otimes E_{21}) + (E_{32}
\otimes E_{23})  + (E_{31} \otimes E_{13}) $
\end{enumerate}

We note that the first three (and similarly the last three) are
connected via odd reflections which correspond to the odd
reflections that connect the associated Dynkin diagrams. The even
reflection which changes the signs of the even roots will
connect the first three to the last three. Hence all these
r-matrices are related to one another via (even or odd)
reflections.

\subsection{Constructing Nonstandard r-matrices}
\label{nonstandardsubsection}
For any given admissible triple ${(\Gamma_1, \Gamma_2,
\tau)}$, we define a partial order on the set of positive roots,
and then according to this setup, the r-matrix we obtain from
our theorem is of the form:  
\begin{equation*}
r = r_0 + 
\sum_{\alpha > 0} e_{-\alpha} \otimes e_{\alpha}  + 
\sum_{\alpha,\beta > 0, \alpha < \beta} 
(e_{-\alpha} \otimes e_{\beta} - (-1)^{|\alpha|} e_{\beta}
\otimes e_{-\alpha})  
\end{equation*}
\noindent 
where the particular $r_0 \in \h \otimes \h$ has to satisfy
${(\tau(\alpha) \otimes 1) (r_0) + (1 \otimes \alpha) (r_0) = 0}$
for all ${\alpha \in \Gamma_1}$. [Of course we still assume
${r_0 + {r_0}^{21} = \Omega_0}$].
\\
\\
\noindent
In our case then, the nonstandard r-matrices come from the two
Dynkin diagrams $D_2$ and $D_5$: 

For $D_2$ let ${\Gamma_1 = \{\alpha_1 + \alpha_2\}}$ and
${\Gamma_2 = \{-\alpha_2\}}$. Define ${\tau(\alpha_1 + \alpha_2)
= -\alpha_2}$. Then the partial order on positive roots will be:
${\alpha_1 + \alpha_2 < -\alpha_2}$. If $r_0$ is given as above
(i.e. ${(-\alpha_2 \otimes 1) (r_0) + (1 \otimes (\alpha_1 +
\alpha_2))(r_0) = 0}$), then the associated r-matrix will be:
\begin{eqnarray*}
r_{ns_1} &=& r_0 + (E_{21} \otimes E_{12}) - (E_{23} \otimes
E_{32})  + (E_{31} \otimes E_{13}) \\
         &+& \left ( (E_{31} \otimes E_{32}) + (E_{32}
\otimes E_{31}) \right )
\end{eqnarray*}
\noindent
The first few terms will actually make up $r_{st}(D_2)$ for the
chosen $r_0$, so we can rewrite the above as:
$$r_{ns_1} = r_{st}(D_2) + (E_{31} \otimes E_{32}) + (E_{32}
\otimes E_{31}) $$

For $D_5$ let ${\Gamma_1 = \{\alpha_2\}}$ and
${\Gamma_2 = \{-\alpha_1 -\alpha_2\}}$. Define ${\tau(\alpha_2)
= -\alpha_1 - \alpha_2}$. Then the partial order on positive
roots will be: ${\alpha_2 < -\alpha_1 -\alpha_2}$. If $r_0$ is
given as above (i.e. ${((-\alpha_1 -\alpha_2) \otimes 1) (r_0) +
(1 \otimes \alpha_2)(r_0) = 0}$), then the associated r-matrix
will be:  
\begin{eqnarray*} 
r_{ns_2} &=& r_0 + (E_{12} \otimes E_{21}) + (E_{32} \otimes
E_{23}) - (E_{13} \otimes E_{31})
\\          
&+& \left ((E_{32} \otimes E_{31}) + (E_{31} \otimes E_{32})
\right ) \end{eqnarray*}
\noindent
The first few terms will actually make up $r_{st}(D_5)$ for the
chosen $r_0$, so we can rewrite the above as:
$$r_{ns_2} = r_{st}(D_5) + (E_{32} \otimes E_{31}) + (E_{31}
\otimes E_{32}) $$

Note that if for $D_2$ we let ${\tau(-\alpha_2) =
\alpha_1 + \alpha_2}$, then we would have: ${-\alpha_2 <
\alpha_1 + \alpha_2}$, and we would get  
$$ r_{ns_3} = r_{st}(D_2) + (-E_{23} \otimes E_{13}) + (-E_{13}\otimes 
E_{23}) ;$$
\noindent 
and if for $D_5$ we let $\tau(-\alpha_1 - \alpha_2) = \alpha_2$,
then the order would be:  ${-\alpha_1 - \alpha_2 < \alpha_2}$,
and we would get 
$$ r_{ns_4} = r_{st}(D_5) + (-E_{13} \otimes E_{23}) + (-E_{23}
\otimes E_{13}). $$

\noindent
Hence the nonstandard r-matrices that can be constructed using our 
theorem are:

\begin{enumerate}
\item $r_{ns_1} = r_{st}(D_2) + (E_{31} \otimes E_{32}) + 
(E_{32} \otimes E_{31});$ 
\item $r_{ns_2} = r_{st}(D_5) + (E_{31} \otimes E_{32}) +
(E_{32} \otimes E_{31});$    
\item $ r_{ns_3} = r_{st}(D_2) + (-E_{13}\otimes E_{23}) +
(-E_{23} \otimes E_{13});$   
\item $ r_{ns_4} = r_{st}(D_5) + (-E_{13} \otimes E_{23}) +
(-E_{23} \otimes E_{13}). $   
\end{enumerate}

\section{Conclusion} 

In the Lie algebra case, the main classification theorem comes in two 
parts. The constructive part that gives an r-matrix for a given 
admissible triple is accompanied with the assertion that any given 
r-matrix that satisfies ${r + r^{21} = \Omega}$ can be obtained by the 
same construction for a suitable choice of an admissible triple. We 
would like to prove such an assertion for Lie superalgebras, or  
come up with a counterexample.

We consider once again the simple Lie superalgebra $sl(2,1)$. We
define:  
$$\begin{matrix} 
f(E_{11}+E_{33}) = 0 & f(E_{22}+E_{33}) = E_{22}+E_{33} \\
f(E_{21}) = 0 & f(E_{12}) = E_{12} \\
f(E_{23}) = 0 & f(E_{13}) = E_{13} \\
f(E_{31}) = -E_{13} & f(E_{32}) = E_{23} + E_{32}  
\end{matrix}$$
\noindent
and extend $f$ to a linear map on $\g$. We can easily
check that this function satisfies  
$$(f-1)[f(x),f(y)] = f([(f-1)(x),(f-1)(y)]),$$
\noindent
which is equivalent to the associated $2-$tensor being an
$r-$matrix. \footnotemark

\footnotetext{Equivalently we can see that the Cayley transform
$\Theta$ is an isomorphism: $\Theta$ maps $\overline{(f-1)(x)}$
to $\overline{f(x)}$. In ${C_1}/{{C_1}^{\perp}}$ the only
nontrivial coset is $\overline{E_{13}+ E_{31}}$, and in
${C_2}/{{C_2}^{\perp}}$ it is $\overline{E_{23}+ E_{32}}$. So we
can choose ${x = E_{32} - E_{13} - E_{31}}$.}

We write the quadratic Casimir element:
$$ \Omega = \Omega_0 + 
( E_{12}\otimes E_{21} + E_{21}\otimes E_{12} ) + 
(-E_{13}\otimes E_{31} + E_{31}\otimes E_{13} ) +
(-E_{23}\otimes E_{32} + E_{32}\otimes E_{23} ) $$
\noindent
where $\Omega_0 = (E_{11}+ E_{33}) \otimes (-E_{22}-E_{33}) + 
 (-E_{22}-E_{33}) \otimes (E_{11}+E_{33}). $
Then if we define $r(f)$ to be the $2-$tensor $(f \otimes
1)\Omega$, we get: 
$$r(f) = r_0 + E_{12}\otimes E_{21} - E_{13} \otimes E_{31} +
E_{32} \otimes E_{23}  - E_{13}\otimes E_{13} +   E_{23} \otimes
E_{23}  $$
\noindent
where $r_0 = (-E_{22} - E_{33}) \otimes (E_{11} + E_{33}).$ It
is clear that $r(f)$ satisfies Equation \ref{unitarity}. 

This r-matrix is not among those constructed using Theorem 
\ref{SuperTheorem}. In fact we can prove that the two subsuperalgebras 
$Im(f)$ and $Im(f-1)$ will never be simultaneously isomorphic to root 
subsuperalgebras. The corresponding subsuperalgebras for functions 
constructed by the theorem will always be root subsuperalgebras. Thus the
Belavin-Drinfeld type data we used is not enough to classify
all solutions to the system of Equations \ref{unitarity} and
\ref{YBequation}. We hope to address this problem in a
separate paper.

\end{document}